\documentclass[aos]{imsart}

\RequirePackage{amsthm,amsmath,amsfonts,amssymb}
\RequirePackage[numbers]{natbib}
\RequirePackage[colorlinks,citecolor=blue,urlcolor=blue]{hyperref}
\RequirePackage{graphicx}

\usepackage{amsmath,amsfonts,amsthm}
\usepackage{mathrsfs}
\usepackage{amssymb, color}

\startlocaldefs
\newtheorem{theorem}{Theorem}[section]

\newtheorem{proposition}[theorem]{Proposition}
\theoremstyle{remark}

\DeclareMathOperator{\E}{\mathbb{E}}
\let\P\relax\DeclareMathOperator{\P}{\mathbb{P}}
\DeclareMathOperator{\Hess}{\mathrm{Hess}}
\DeclareMathOperator{\Var}{\mathrm{Var}}

\def\calS{{\cal S}}

\def\calA{{\cal A}}
\def\calC{{\cal C}}
\def\calL{{\cal L}}

\def\calH{{\cal H}}

\def\calG{{\cal G}}

\def\calN{{\cal N}}

\def\R{{\mathbb R}}

\def\E{{\mathbb E}}

\def\scrF{{\mathscr F}}

\endlocaldefs

\begin{document}

\begin{frontmatter}
\title{Stein's method of normal approximation: \\ Some recollections and reflections}
%\title{A sample article title with some additional note\thanksref{t1}}
\runtitle{Stein's method of normal approximation}
%\thankstext{T1}{A sample additional note to the title.}

\begin{aug}
%%%%%%%%%%%%%%%%%%%%%%%%%%%%%%%%%%%%%%%%%%%%%%
%%Only one address is permitted per author. %%
%%Only division, organization and e-mail is %%
%%included in the address.                  %%
%%Additional information can be included in %%
%%the Acknowledgments section if necessary. %%
%%%%%%%%%%%%%%%%%%%%%%%%%%%%%%%%%%%%%%%%%%%%%%
\author[A]{\fnms{Louis H.\,Y.} \snm{Chen}\ead[label=e1]{matchyl@nus.edu.sg}},
%\author[B]{\fnms{Second} \snm{Author}\ead[label=e2,mark]{second@somewhere.com}}
%\and
%\author[B]{\fnms{Third} \snm{Author}\ead[label=e3,mark]{third@somewhere.com}}
%%%%%%%%%%%%%%%%%%%%%%%%%%%%%%%%%%%%%%%%%%%%%%
%% Addresses                                %%
%%%%%%%%%%%%%%%%%%%%%%%%%%%%%%%%%%%%%%%%%%%%%%
\address[A]{Department of Mathematics,
National University of Singapore,
\printead{e1}}

%\address[B]{Department,
%University or Company Name,
%\printead{e2,e3}}
\end{aug}

\begin{abstract} This paper is a short exposition of Stein's method of normal approximation from my personal perspective. It focuses mainly on the characterization of the normal distribution and the construction of Stein identities. Through examples, it provides glimpses into the many approaches to constructing Stein identities and the diverse applications of Stein’s method to mathematical problems. It also includes anecdotes of historical interest, including how Stein discovered his method and how I found an unpublished proof of his of the Berry-Esseen theorem.
\end{abstract}



\end{frontmatter}
%%%%%%%%%%%%%%%%%%%%%%%%%%%%%%%%%%%%%%%%%%%%%%
%% Please use \tableofcontents for articles %%
%% with 50 pages and more                   %%
%%%%%%%%%%%%%%%%%%%%%%%%%%%%%%%%%%%%%%%%%%%%%%
%\tableofcontents

\setcounter{footnote}{1}
\begin{center}
\begin{minipage}{0.45\textwidth}
\small\it
A thinker original and independent,\\
In search of perfection invariant,\\
Found admissible wisdom's counterexample,\\
Made (fame, humility) exchangeable.\footnotemark
\end{minipage}
\end{center}
\footnotetext{Written by Louis Chen and Yu-Kiang Leong, presented to Charles Stein on his 90th birthday.}

\section{Introduction} Charles Stein was one of the most original statisticians and probabilists of the 20th century. His ideas and his work have profoundly influenced the development of both statistics and probability. Perhaps his greatest contribution to probability is his novel method of proving central limit theorems with explicit finite sample error bounds, which we now call Stein's method. Stein started to develop this new method in the 1960s, when he was trying to prove the combinatorial central limit theorem using his own approach. But it was not until 1972 that his groundbreaking paper "A bound for the error in the normal approximation to the distribution of a sum of dependent random variables" \cite{Stein1972} was published. In this paper, an abstract normal approximation theorem is proved, which is then applied to obtain Berry-Esseen bounds for the accuracy of the normal approximation to the distribution of the sum of a stationary sequence of random variables. 

A striking feature of Stein’s method of proof is that it does not make any use of Fourier analytic methods, but works directly on spaces of random variables with the solution of a first order differential equation playing a critical role. The flexibility of the method is also evident, since the dependence setting of stationary sequences is quite different from that of the combinatorial central limit theorem, with which Stein began.
 
Stein expanded his method into a definitive theory in his monograph, {\it Approximate Computation of Expectations} \cite{Stein1986}, where he developed the notion of an exchangeable pair of random variables and pointed the way to approximation by distributions other than the normal. Indeed, Stein's ideas are very general and, in principle, can be applied to any distributional approximation. Since the publication of his 1972 paper, many approximations by other distributions have been developed, starting from the classical ones, such as the Poisson (Chen \cite{Chen1975}), the bionomial (Stein \cite{Stein1986}) and the gamma (Luk \cite{Luk1994}) distributions, and moving on to the exotic and lesser known cases, such as the Dickman distribution (Bhattacharjee and Goldstein \cite{B-G2019}), and the list continues to grow.  
 
In this paper, I will give a short exposition of Stein's method of normal approximation from my personal perspective. I will focus mainly on the characterization of the normal distribution and the construction of Stein identities. Through examples, I will provide glimpses into the many approaches to constructing Stein identities and the diverse applications of Stein’s method to mathematical problems even within normal approximation. I will also include some personal anecdotes of historical interest, including how Stein discovered his method and how I found an unpublished proof of his of the Berry-Esseen theorem.

\section{From characterization to approximation} \label{characterization}  A way to understand Stein's method of normal approximation is to begin with Stein's characterization of the normal distribution, which states that for a random variable $W$ to have the standard normal distribution, it is necessary and suffcient that 
\begin{equation} \label{Stein-id-normal}
\E\{f'(W) - Wf(W)\} = 0 ~~ \text{for}~~ f \in \mathcal{G},
\end{equation}
where $\mathcal{G}$ is the class of all absolutely continuous functions $f: \R \rightarrow \R$ with a.e. derivative $f'$ such that $\E| f'(Z)| < \infty$, and $Z \sim \mathcal{N}(0,1)$. A proof of this result can be found on p. 21 in Stein's 1986 monograph \cite{Stein1986}. This characterization of the normal distribution implies that if a random variable $W$ with $\E W = 0$ and $\mathrm{Var}(W) = 1$ is not distributed as $\mathcal{N}(0,1)$, then $\E\{f'(W) - Wf(W)\} \neq 0$ for some $f \in \mathcal{G}$. The discrepancy between $\mathcal{L}(W)$ and $\mathcal{N}(0,1)$ can then be measured by letting $f = f_h$, a bounded solution of the first order differential equation, 
\begin{equation} \label{Stein-eqn}
f'(w) - wf(w) = h - \E h(Z),
\end{equation}
where the test function $h \in \mathcal{H}$, a separating class of functions such that $f_h \in \mathcal{G}$. By a separating class $\mathcal{H}$, we mean a class of Borel measurable functions $h: \R \rightarrow \R$ such that two random variables, $X$ and $Y,$  have the same distribution if $\E h(X) = \E h(Y)$ for $h \in \mathcal{H}$. Such a separating class contains functions $h$ for which both $\E h(X)$ and $\E h(Y)$ exist.

The equation (\ref{Stein-eqn}), which we call the Stein equation, then gives 
\begin{equation} \label{Eh-Eh}
\E h(W) - \E h(Z) = \E\{f_h'(W) - Wf_h(W)\}.
\end{equation}
Define a distance between $\calL(W)$ and $\calN(0,1)$ induced by $\calH$ as 
\begin{equation} \label{def-d-calH}
d_\calH(W,Z) := \sup_{h\in\calH}|\E h(W) - \E h(Z)| = \sup_{h\in\calH}|\E\{f_h'(W) - Wf_h(W)\}|.
\end{equation}
Then bounding $d_\calH(W,Z)$ is equivalent to bounding $\sup_{h\in\calH}|\E\{f_h'(W) - Wf_h(W)\}|$, for which one uses the probabilistic structure of $W$ and the boundedness and smoothness properties of $f_h$. The remarkable power of Stein's method is that $\sup_{h\in\calH}|\E\{f_h'(W) - Wf_h(W)\}|$ is a lot easier to handle than $\sup_{h\in\calH}|\E h(W) - \E h(Z)|$. If $\calL(W)$ is close to $\calN(0,1)$, we should expect $\sup_{h\in\calH}|\E\{f_h'(W) - Wf_h(W)\}|$ to be small.

Three seperating classes of functions are of interest. These are:
\begin{eqnarray*}
\calH_W &=& \{h: |h(u)-h(v)| \le |u-v|\}, \\
\calH_K &=& \{h: h(w) =1 ~\text{for} ~w\le x ~ \text{and} = 0  ~\text{for} ~ w> x, x \in \R\}, \\
\calH_{TV} &=& \{h: h(w) = I(w\in A), ~A ~ \text{is a Borel subset of}~ \R\}.
\end{eqnarray*}
We define the Wasserstein distance, denoted by $d_W$, to be the distance induced by $\calH_W$, the Kolmogorov distance, denoted by $d_K$, to be the distance induced by $\calH_K$, and the total variation distance, denoted by $d_{TV}$, to be the distance induced by $\calH_{TV}$. Note that $d_K(W,Z) = \sup_{x \in \R}|\P(W \le x) - \Phi(x)|$, where, and for the rest of this paper, $\Phi$ denotes the standard normal distribution function.

If we define the operator $L$ by $Lf(w) = f'(w) - wf(w)$, then (\ref{Stein-id-normal}) can be written as $\E Lf(W) =0 ~\text{for}~ f \in \calG$. We call $L$ a Stein operator for $\calN(0,1)$ and say $L$ characterizes $\calN(0,1)$.

For the rest of this paper, $Z \sim \calN(0,1)$. Unless differently specified, all functions are from $\R$ to $\R$. The following are the classes of functions used: (i)  $\calC = \{f: f~ \text{is continuous}\}$, (ii) $\calC^k = \{f: f, f', \cdots, f^{(k)} ~ \text{are continuous}\}$, (iii) $\calC_B = \{f: f ~ \text{is bounded and continuous}\}$, (iv) $\calC_B^k = \{f: f, f',\cdots, f^{(k)} ~ \text{are bounded and continuous}\}$.

The following two propositions provide alternative definitions of the Wasserstein and total variation distances and can be very useful in applications. Propostion \ref{prop-TV} is proved using Lusin's theorem (see, for example, Rudin \cite{Rudin1986}, p. 53). 
\begin{proposition} \label{prop-W}
We have
\begin{equation}
d_W(W,Z)= \sup_{h \in \mathcal{C}^1, \|h'\|_\infty \le 1} |\E h(W)-\E h(Z)|.
\end{equation}
\end{proposition}

\begin{proposition} \label{prop-TV}
We have
\begin{equation} \label{d-TV-h-B}
d_{TV}(W,Z) =  \sup_{0\le h \le1} |\E h(W)-\E h(Z)| = \sup_{h\in \mathcal{C}, 0\le h \le 1}|\E h(W)-\E h(Z)|.
\end{equation}
\end{proposition}

\section{Stein identities and error terms} \label{Stein-id}
How do we use the probability structure of $W$ for bounding the right hand side of (\ref{def-d-calH}), namely, $\sup_{h\in\calH}|\E\{f_h'(W) - Wf_h(W)\}|$? A way to do this is by constructing what we call a Stein identity. This is perhaps best understood by looking at a specific example.

Let $X_1, \cdots,X_n$ be independent random variables with $\E X_i = 0$, $\mathrm{Var}(X_i) =\sigma_i^2$ and $\E|X_i|^3 < \infty$. Let $W = \sum_{i=1}^n X_i$ and $W^{(i)} = W - X_i$ for $i=1,\cdots,n$. Assume that $\mathrm{Var}(W) = 1$, which implies $\sum_{i=1}^n \sigma_i^2 = 1$. Let $f$ be absolutely continuous such that both $f$ and $f'$ are bounded. Using the independence among the $X_i$ and the property that $\E X_i = 0$, and also Fubini's theorem, we have
\begin{eqnarray} \label{Stein-identity-indep-1}
\E Wf(W) &=& \sum_{i=1}^n \E X_i f(W) =\sum_{i=1}^n \E X_i [f(W^{(i)}+X_i) - f(W^{(i)})] \\ \notag
&=& \sum_{i=1}^n\E\int_0^{X_i} X_i f'(W^{(i)}+t) dt = \sum_{i=1}^n\E\int_{-\infty}^\infty f'(W^{(i)}+t)\hat{K}_i(t)dt \\  \notag
&=& \sum_{i=1}^n\E\int_{-\infty}^\infty f'(W^{(i)}+t)K_i(t)dt
\end{eqnarray}
where 
\begin{equation*}
\hat{K}_i(t) = X_i[I(X_i>t>0)-I(X_i<t\le 0)], ~  K_i(t) = \E\hat{K}_i(t).
\end{equation*}
It can be shown that for each $i$, $\sigma_i^{-2}K_i$ is a probability density function.  So (\ref{Stein-identity-indep-1}) can be rewritten as 
\begin{equation} \label{Stein-identity-indep-2}
\E Wf(W) = \sum_{i=1}^n \sigma_i^2\E f'(W^{(i)}+T_i)
\end{equation}
where $T_1,\cdots,T_n, X_1,\cdots,X_n$ are independent and $T_i$ has the density $\sigma_i^{-2}K_i$.
Both the equations (\ref{Stein-identity-indep-1}) and (\ref{Stein-identity-indep-2}) are Stein identities for $W$. (The equation $\E Zf(Z) = \E f'(Z)$ for $f \in \calG$ is a Stein equation for $\mathcal{N}(0,1)$.)
From (\ref{Stein-identity-indep-2}) we obtain
\begin{eqnarray} \label{Stein-identity-indep-3}
\E [f'(W) - Wf(W)]&=& \sum_{i=1}^n \sigma_i^2\E[f'(W) - f'(W^{(i)})] \\  \notag
&& - \sum_{i=1}^n \sigma_i^2\E[f'(W^{(i)}+T_i) - f'(W^{(i)})] 
\end{eqnarray}
where the error terms on the right hand side provide an expression for the deviation of $\E[f'(W) - Wf(W)]$ from $0$.
Now let $f$ be $f_h$ where the test function $h \in \mathcal{C}^1$ is such that $\|h'\|_\infty \le 1$. Then $f_h \in \mathcal{C}_B^2$. We show how easily (\ref{Stein-identity-indep-3}) leads to a bound on the Wasserstein distance between $\calL(W)$ and $\calN(0,1)$, by applying Taylor's theorem to the right hand side of (\ref{Stein-identity-indep-3}), and using Proposition \ref{prop-W} and the fact that $\|f''_h\|_{\infty} \le 2\|h'\|_{\infty}$ (see Lemma 2.4 in Chen, Goldstein and Shao \cite{C-G-S2011}):
\begin{eqnarray*} 
d_W(W,Z) &\le& \sum_{i=1}^n \sup_{h\in \mathcal{C}^1, \|h'\|_\infty \le 1} \sigma_i^2\E|X_i| \|f_h''\|_\infty    
+\sum_{i=1}^n \sup_{h\in \mathcal{C}^1, \|h'\|_\infty \le 1} \sigma_i^2\E|T_i|\|f_h''\|_\infty \\ 
&\le& 2\sum_{i=1}^n\sigma_i^2\E|X_i| +2\sum_{i=1}^n\sigma_i^2\E|T_i|.
\end{eqnarray*}
Since $\sigma_i^2\E|X_i| \le (\E|X_i|^3)^{2/3}(\E|X_i|^3)^{1/3} = \E|X_i|^3$ and $\sigma_i^2\E|T_i| = \sigma_i^2(1/2\sigma_i^2)\E|X_i|^3 = (1/2)\E|X_i|^3$, we have
\begin{equation} \label{bd-d-W-indep}
d_W(W,Z) \le 3\sum_{i=1}^n\E|X_i|^3.
\end{equation}
The order of the bound in (\ref{bd-d-W-indep}) is optimal. 

If we let the test function $h$ be the indicator function $h_x$ given by $h_x(w) =I(w \le x)$ where $x \in \R$, and denote  $f_{h_x}$ by $f_x$, we obtain from (\ref{Stein-identity-indep-2}) the Kolmogorov distance between $\calL(W)$ and $\calN(0,1)$ as given below
$$
d_K(W,Z) = \sup_{x \in \R}\sum_{i=1}^n \sigma_i^2\E[f_x'(W) - f_x'(W^{(i)} + T_i)].
$$
As $f_x'$ is discontinuous at $x$, we cannot apply Taylor's theorem. We need to develop different techniques to bound the Kolmogorov distance. One such technique is the concentration inequality approach developed by Stein, which we will discuss in section \ref{concentation-ineq}.

The above method of constructing Stein identities has been developed for locally dependent random variables (see, for example, Chen and Shao \cite{Chen-Shao2004}).

 In \cite{Stein1972}, Stein used auxilliary randomization to construct Stein identities, and in \cite{Stein1986, Stein1992}, he introduced the notion of exchangeable pair and used auxilliary randomization to construct exchangeable pairs of random variables, which he in turn used to construct Stein identities. I will discuss exchangeable pairs in section \ref{ex-pair}.

Goldstein and Reinert \cite{Goldstein-Reinert1997} introduced the notion of the zero-bias transformation where any random variable $W$ with $\E W =0$ and $\mathrm{Var}(W) = B^2 > 0$ has a zero-bias transform $W^*$ such that 
\begin{equation} \label{Stein-id-zb}
\E Wf(W) = B^2\E f'(W^*)
\end{equation}
for absolutely continuous $f$ such that $\E|f'(W)| < \infty$. The Stein identity  (\ref{Stein-id-zb}) is proved by applying Fubini's theorem, and $W^*$ is absolutely continuous with density function given by $\E WI(W > x)/B^2$. If $W$ and $W^*$ are defined on the same probability space and $B^2 =1$, we can use (\ref{Stein-id-zb}) to obtain quite effortlessly
\begin{equation}\label{bd-d-W-zb}
d_W(W,Z) \le 2\E|W^*-W|.
\end{equation}
The construction of $W^*$ can be tricky. Nevertheless, in the case of  independent $X_1,\cdots,X_n$ considered above, 
$W^* = W^{(I)}+X_I^*$ where $I$, $X_i$, ${X_i^*}$, $i =1,\cdots, n$, are independent and $\P(I = i) = \sigma^2$. In fact in (\ref{Stein-identity-indep-2}), $T_i = X_i^*$. So applying (\ref{bd-d-W-zb}), we obtain (\ref{bd-d-W-indep}). Also, in \cite{Goldstein2010}, Goldstein made clever use of zero-bias coupling to improve the constant in the bound in (\ref{bd-d-W-indep}) from $3$ to $1$.

Closely related to zero-bias transformation is the size-bias transformation of non-negative random variables. Every nonnegative random variable $Y$ with mean $\mu$ has a size-bias transform $Y^s$ such that  
\begin{equation} \label{Stein-id-sb1}
\E Yf(Y) = \mu\E f(Y^s)
\end{equation}
for $f$ for which the expectations exist. Assuming that $\Var(Y) = B^2$ and that $Y$ and $Y^s$ are defined on the same probability space, we can define $W = (Y - \mu)/B$ and $W^s = (Y^s - \mu)/B$, and construct a Stein identity for $W$ as follows: for absolutely continuous $f$ such that $f$ and $f'$ are bounded,
\begin{eqnarray} \label{Stein-id-sb2}
\E Wf(W) &=& \frac{\mu}{B}\E \{f(W^s) - f(W)\} = \frac{\mu}{B}\E\int_{0}^{W^s-W}f'(W+t)dt \\ \notag
&=& \E\int_{-\infty}^\infty f'(W+t)\hat{K}(t)dt,
\end{eqnarray}
where
$$
\hat{K}(t) = \frac{\mu}{B}[I(W^s-W > t >0) - I(W^s-W < t \le 0)]. 
$$
The method of size-bias coupling in Stein's method was introduced in Baldi, Rinott and Stein \cite{B-R-S1989} and fully developed in Goldstein and Rinott \cite{Goldstein-Rinott1996}.

For a random measure $\Xi$ with finite intensity measure $\Lambda$ on a locally compact separable metric space $\Gamma$, the  Palm measure $\Xi_\alpha$ associated with $\Xi$ at $\alpha \in \Gamma$ (see Kallenberg \cite{Kallenberg1983}, pp. 83, 103) can be thought of as a size-bias transform of $\Xi$. If $\Xi$ is a simple point process, the distribution of $\Xi_\alpha$ can be interpreted as the conditional distribution of $\Xi$ given that a point of $\Xi$ at $\alpha$ has occurred. On the other hand, if $\Lambda(\{\alpha\})>0$, then $\Xi(\{\alpha\})$ is a non-negative random variable with positive mean and $\Xi_\alpha(\{\alpha\})$ is the size-bias transform of $\Xi(\{\alpha\})$ in the sense of (\ref{Stein-id-sb1}).

From Palm theory (see Kallenberg \cite{Kallenberg1983}, p. 84), we have for absolutely continuous $f: \R \rightarrow \R$ such that the expectations and the integral in (\ref{Palm-|Xi|}) exist,

\begin{equation} \label{Palm-|Xi|}
\E|\Xi|f(|\Xi|)= \E\int_\Gamma f(|\Xi_\alpha|)\Lambda(d\alpha),
\end{equation}
where $|\Xi|=\Xi(\Gamma)$. Assume that $\Xi$ and $\Xi_\alpha$, $\alpha \in \Gamma$, are defined on the same probability space. Let $I$ be a random element taking values in $\Gamma$ with probability measure $\Lambda/\lambda$, where $\lambda = \Lambda(\Gamma)$, and be independent of $\Xi$ and $\Xi_\alpha$, $\alpha \in \Gamma$. Then (\ref{Palm-|Xi|}) can be written as $\E|\Xi|f(|\Xi|) = \lambda\E f(|\Xi|_I)$, so that $|\Xi_I|$ is the size bias transform of $|\Xi|$ in the sense of (\ref{Stein-id-sb1}), and (\ref{Stein-id-sb2}) is a Stein identity for $W = (|\Xi|-\lambda)/B$ with $W^s = (|\Xi_I| - \lambda)/B$ and $B^2 = \Var(W) < \infty$.

Chen, R\"{o}llin and Xia \cite{C-R-X2021} considered a Stein identity of the form
\begin{equation} \label{Stein-id-random-K}
\E Wf(W) = \E \int_\R f'(W + t)\hat{K}(t)dt,
\end{equation}
where $\hat{K}$ is a random function and obtained a general bound on the Kolmogorov distance in the normal approximation for $W$. They then applied the result to random measures with application to stochastic geometry, and also to the Stein couplings of Chen and R\"{o}llin \cite{Chen-Rollin2010}. Note that the Stein identity (\ref{Stein-id-random-K}) is more general than (\ref{Stein-id-sb2}) as there is no assumption on the $\hat{K}$ in (\ref{Stein-id-random-K}) other than it being a random function.

Chatterjee \cite{Chatterjee2009, Chatterjee2014} considered $W = g(X)$ with $\E W = 0$ and $\mathrm{Var}(W) = 1$, where $X = (X_1,\cdots,X_n)$ is a vector of independent standard normal random variables, and $g: \R^n \rightarrow \R$ is a twice continuously differentiable function with gradient $\nabla g$ and Hessian matrix $\Hess g$. Using Gaussian interpolation, he constructed this Stein identity
\begin{equation} \label{Stein-id-g}
\E Wf(W) = \E T(X)f'(W),
\end{equation}
where $f$ is absolutely continuous such that $\E|Wf(W)| < \infty$ and $f'$ is bounded, and 
$$
T(x) = \int_0^1\frac{1}{2\sqrt{t}}\E\left(\sum_{i=1}^n\frac{\partial g}{\partial x_i}(x)\frac{\partial g}{\partial x_i}(\sqrt{t}x + \sqrt{1-t}X\right)dt.
$$
From (\ref{Stein-id-g}), by letting $f(w) = w$, we obtain  $\E T(X) = 1$ .  Let $h$ be continuous such that $0 \le h \le 1$. Then by Lemma 2.4 in Chen, Goldstein and Shao \cite{C-G-S2011}, $\|f'_h\|_{\infty} \le 2\|h\|_{\infty} \le 2$. From (\ref{Stein-id-g}) and applying Proposition \ref{prop-TV}, we obtain
\begin{equation} \label{d-TV-g}
d_{TV}(W, Z) = \sup_{h \in \calC, 0\le h \le1}|\E(1 - T(X))f_h'(W)| \le 2\E|T(X) -1| \le 2\sqrt{\Var(T(X))}.
\end{equation}
where for the last inequality $T(X)$ is assumed to be square integrable.

To bound $\sqrt{\Var(T(X))}$, Chatterjee \cite{Chatterjee2009} developed a Poincar\'e inequality on $\Var(T(X))$, where the upper bound involves the Euclidean norm of $\nabla g$ and the operator norm of $\Hess g$. He called this Poincar\'e inequality a \textit{second order Poincar\'e inequality}. This inequality was used to prove central limit theorems for linear statistics of eigenvalues of random matrices.

Nourdin and Peccati \cite{N-P2009} established a fundamental connection between Stein's method of normal approximation and Malliavin calculus. Let $B = (B_t)_{t \in \R_+}$ be a standard Brownian motion on a complete probability space $(\Omega,\scrF, P)$, where $\scrF$ is generated by  $(B_t)_{t \in \R_+}$. Let $F$ be a functional of Gasussian processes lying in a suitable subspace of $L^2(\Omega) = L^2(\Omega,\scrF, P)$ such that $\E F = 0$ and $\mathrm{Var}(F) = 1$. Let $L$ be the Ornstein-Uhlenback operator with pseudo-inverse $L^{-1}$ and let $f \in \calC^1$ such that $\E|Ff(F)| < \infty$ and $f'$ is bounded. Using integration by parts, Nourdin and Peccati \cite{N-P2009} obtained the following Stein identity,
\begin{equation} \label{Stein-id-F}
\E F f(F) = \E [f'(F)\langle DF,-DL^{-1}F\rangle_{L^2(\R_+)}] = \E Tf'(F),
\end{equation}
where $\langle \cdot, \cdot \rangle_{L^2(\R_+)}$ is the inner product on $L^2(\R_+, dx)$, $D$ the Malliavin derivative and $T = \langle DF,-DL^{-1}F\rangle_{L^2(\R_+)}$. The $F$ considered in (\ref{Stein-id-F}) generalizes the result of Chatterjee \cite{Chatterjee2009} to settings beyond the Euclidean space.

From (\ref{Stein-id-F}), $d_{TV}(F,Z) \le 2\sqrt{\Var(T)}$, where $T$ is assumed to be square integrable. Nourdin and Peccati \cite{N-P2009} applied this result to a few specific functionals including proving the fourth moment theorem which provides total variation error bounds of the order $\sqrt{\E(F_n^4) - 1}$ in the central limit theorem for elements $\{F_n\}_{n \ge 1}$ of Wiener chaos of any fixed order such that $\E F_n^2 = 1$.

All the Stein identities discussed above are proved using the principle of integration by parts. This is not surprising since the 
characterization equation (\ref{Stein-id-normal}) is proved using Fubini's theorem, which is more general than integration by parts, if $W \sim \calN(0,1)$. 

\section{From teaching, a creative spark} I do not recall having asked Stein about how he discovered his new method of normal approximation when I was his Ph.D. student at Stanford University. It was years later through conversations with him and Persi Diaconis that I learned that his method grew out of his teaching a course of nonparametric statistics in the 1960s. He was proving the combinatorial central limit theorem in the course and decided to try his own approach instead of presenting the published work of Wald and Wolfowitz \cite{W-W1944} or of Hoeffding \cite{Hoeffding1951}. 

Let $W_n = \sum_{i=1}^n c_{n,i\pi(i)}$ with $\E W_n =0$ and $\mathrm{Var}(W_n)= 1$, where $\pi$ is a random permutation of $\{1,2,\dots,n\}$ and $\{c_{n,ij}: i,j = 1, \dots, n\}$ is an $n \times n$  array of real numbers such that for every $i$, $\sum_{j=1}^nc_{n,ij} =0$ and for every $j$, $\sum_{i=1}^nc_{n,ij} = 0$. Let $\phi_n$ be the characteristic function of $W_n$. Using the notion of exchangeable pair, Stein showed that for $t \in \R$, 
\begin{equation} \label{phi-n}
\phi_n'(t) + t\phi_n(t) \rightarrow 0 ~ \text{as} ~ n \rightarrow \infty.
\end{equation}
Then the characteristic function $\phi$ of the weak limit of any subsequence of $W_n$ satisfies 
\begin{equation} \label{phi}
\phi'(t) + t\phi(t) = 0, 
\end{equation}
which implies $\phi(t) = e^{-t^2/2}$. This in turn implies $W_n \stackrel{\mathcal{L}}{\rightarrow} \mathcal{N}(0,1)$.

Stein later realized that there was nothing special about the complex exponential in the characteristic function. Since 
\begin{eqnarray*}
\phi_n'(t) &=& \frac{\partial}{\partial t}\E e^{itW_n} = \E\frac{\partial}{\partial t}e^{itW_n} = i\E W_ne^{itW_n}, \\
t\phi_n(t) &=& t\E e^{itW_n} =  \frac{1}{i}\E\frac{\partial}{\partial W_n}e^{itW_n},
\end{eqnarray*}
replacing the complex exponential $e^{it\cdot}$ by an arbitrary function $f(\cdot)$, (\ref{phi-n}) becomes
$$
\E\{f'(W_n) - W_nf(W_n)\} \rightarrow 0 ~ \text{for a suitable class of}~ f ~\text{as} ~ n \rightarrow \infty.
$$
And the implication of $\phi(t) = e^{-t^2/2}$ from (\ref{phi}) becomes
$$
\E\{f'(Z) - Zf(Z)\} = 0 ~ \text{for a suitable class of}~ f~\text{characterizes} ~ Z \sim \mathcal{N}(0,1).
$$
This then led to calculating $\E\{f'_h(W_n) - W_nf_h(W_n)\}$ for assessing the rate of convergence as described in Section \ref{Stein-id}.

\section{The concentration inequality approach} \label{concentation-ineq} When I became Stein's Ph.D. student in 1969, Stein asked me to read the book {\it Topological Groups} by Leopoldo Nachbin. It seemed that he wanted me to work on decision theory as I had taken his course on decision theory but had not taken any course on Stein's method. At that time, a fellow Ph.D. student, Richard Shorrock, who was my office mate and had become a good friend of mine, encouraged me to work on Stein's method. He passed me a set of lecture notes, which he said were given to him by Roberto Mariano, who took the notes in a course on Stein's method of normal approximation conducted by Stein himself.

I started to read the notes and found that many lectures were incomplete and there were a few incomplete proofs of the Berry-Esseen theorem using Stein's method. This reminded me of my experience when I took Stein's course on decision theory. He would come to class to teach the latest topic in his research but would often get stuck halfway through his proof or calculation. At the end of his lecture, he would say that he would continue his proof or calculation in the next lecture. But when he returned, he would abandon the topic or the approach and start a new one. This method of teaching may be unconventional but it offers challenging problems and research opportunities for students.

I tried to complete the proofs of the Berry-Esseen theorem in the notes but without success. But as I read on, I discovered to my surprise a complete proof of the Berry-Esseen theorem for independent and identically distributed random variables, in which Stein proved what he called a \textit{concentration inequality} and used it to overcome the difficulty arising from the discontinuity of the test function $h$. As this proof has been presented in Ho and Chen \cite{Ho-Chen}, I will give a sketch of it here.

Let $\xi_1, \cdots, \xi_n$ be independent and identically distributed random variables with zero mean, variance $\sigma^2$ and absolute third moment $\gamma$.  Let $W_n = (\sum_{i=1}^n\xi_i)/\sigma\sqrt{n}$. The Berry-Esseen theorem states that there is a constant $C$ such that
\begin{equation} \label{B-E-theorem}
\sup_{x\in \R}|\P(W_n\le x) - \Phi(x)| \le \frac{C\gamma}{\sigma^3\sqrt{n}}.
\end{equation}

Let $X_i = \xi_i/\sigma\sqrt{n}$, $i = 1, \cdots, n$. Then $X_1, \cdots, X_n$ are independent and identically distributed with $\E X_i =0$, $\mathrm{Var}(X_i)=1/n$ and $\E|X_i|^3 = \gamma/\sigma^3 n^{3/2}$. Let $W_k = \sum_{i=1}^k X_i$. So $\mathrm{Var}(W_n) = 1$. From (\ref{Stein-identity-indep-1}), we have, for absolutely continuous and bounded $f$ with bounded $f'$, the Stein identity 
\begin{equation} \label{Stein-id-K}
\E W_nf(W_n) = \E\int_{-\infty}^\infty f'(W_{n-1} +t)K(t)dt
\end{equation}
where $K(t) = n\E X_n[I(X_n > t > 0) - I(X_n \le t \le 0)]$.

As $K$ is a probability density, let $T \sim K$ and be independent of $X_1, \cdots, X_n$. Then $\E|T| = n\E|X_1|^3/2 = \gamma/2\sigma^3\sqrt{n}$ and the Stein identity (\ref{Stein-id-K}) can be written as 
\begin{equation} \label{Stein-id-T}
\E W_nf(W_n) = \E f'(W_{n-1}+T).
\end{equation}
Let $f_x$ be the unique bounded solution of the equation $f'(w) - wf(w) = I(w\le x) - \Phi(x)$, where $x \in \R$. Note that the test fucntion $h$ in this case is an indicator function given by $h(w) = I(w \le x)$. Then, 
\begin{eqnarray} \label{Stein-id-iid}
&&\P(W_n \le x) - \Phi(x) = \E[f_x'(W_n) - W_nf_x(W_n)] ~~~~~~~~~~~~~~~~\\ \notag
&&~~~~=\E[f_x'(W_n) - f_x'(W_{n-1}+T)] \\  \notag
&&~~~~= \P(x-T <  W_{n-1} \le x -X_n, X_n\le T) \\  \notag
&&~~~~~~~ - \P(x-X_n <  W_{n-1} \le x -T, X_n > T) \\  \notag
&&~~~~~~~ + \E[W_nf_x(W_n) - (W_{n-1}+T)f_x(W_{n-1}+T)].
\end{eqnarray}
Noting that
\begin{eqnarray*}
&&|\P(x-T <  W_{n-1} \le x -X_n, X_n\le T) -  \P(x-X_n <  W_{n-1} \le x -T, X_n > T)| \\
&& ~~~~ \le \E \P(x - \max(X_n,T) \le W_{n-1} \le x - \min(X_n,T)|X_n, T),
\end{eqnarray*}
and that 
\begin{eqnarray*}
&&|\E[W_nf_x(W_n) - (W_{n-1}+T)f_x(W_{n-1}+T)]| \\
&& ~~~~\le |\E[W_{n-1}(f_x(W_n) - f_x(W_{n-1}+T))]| +|\E X_nf_x(W_n)| + |\E Tf(W_{n-1} + T)| \\
&& ~~~~ \le \|f_x'\|_\infty \E|W_{n-1}||X_n - T| + \|f_x\|_\infty \E|X_n| + \|f_x\|_\infty \E|T| \\
&& ~~~~ \le (\|f_x'\|_\infty + \|f_x\|_\infty)(\E|X_n| + \E|T|),
\end{eqnarray*}
and using $\|f_x'\|_\infty \le 1$ and $\|f_x\|_\infty \le 1$ (see Lemma 2.3 in Chen, Goldstein and Shao \cite{C-G-S2011}), we obtain from (\ref{Stein-id-iid}),
\begin{eqnarray} \label{P-Phi}
|\P(W_n \le x) - \Phi(x)| &\le&\E \P(x - \max(X_n,T) \le W_{n-1} \le x - \min(X_n,T)|X_n, T) \\ \notag
&& + 3n\E|X_1|^3.
\end{eqnarray}

Let $a < b$ and $\delta > 0$. Let $f$ be such that $f'(w) = I(a - \delta \le w \le b + \delta)$ and $f((a+b)/2) = 0$. Then ${\displaystyle |f| \le \frac{b-a+2\delta}{2}}$. Putting this $f$ in the Stein indentity (\ref{Stein-id-T}) gives

\begin{eqnarray*}
\frac{b-a+2\delta}{2} &\ge& \E W_nf(W_n) = \P(a-\delta \le W_{n-1}+T \le b+\delta) \\
&\ge& \P(a \le W_{n-1} \le b, -\delta \le T \le \delta) \\
&=& \P(a \le W_{n-1} \le b)P(|T| \le \delta) \\
&=& \P(a \le W_{n-1} \le b)(1 - P(|T| > \delta)) \\
&\ge& \P(a \le W_{n-1} \le b)(1 - \E|T|/\delta).
\end{eqnarray*}
Letting $\delta = 2\E|T|$, we obtain this concentration inequality,
\begin{eqnarray} \label{concentration-ineq}
\P(a \le W_{n-1} \le b) &\le& b - a + 4\E|T| = b - a + 2n\E|X_1|^3.
\end{eqnarray}
Combining (\ref{P-Phi}) and (\ref{concentration-ineq}) yields
\begin{equation} \label{BE-bound1}
\sup_{x\in \R}|\P(W_n \le x) - \Phi(x)| \le  6.5n\E|X_1|^3 = \frac{6.5\gamma}{\sigma^3\sqrt{n}}.
\end{equation}

Stein's method would be somewhat incomplete if one could not produce a bound on the Kolmogorov distance depending on the third moment. The concentration inequality approach is a way to achieve that end. This specific piece work of Stein has influenced my research direction for many years. See, for example, \cite{Chen1986, Chen1998, Chen-Fang2015-1, Chen-Fang2015-2, Chen-Shao2001, Chen-Shao2004, Chen-Shao2007}.

\section{Exchangeable pair} \label{ex-pair} Stein had the notion of exchangeable pair and knew how to use it when he was proving the combinatorial central limit therem in the 1960s. In Diaconis \cite{Diaconis1977}, an exchangeable pair of random variables was used in the proof of a normal approximation result for the number of ones in the binary expansion of a random integer. Diaconis wrote that it was a joint work with Stein. It was not until 1986 that Stein in his monograph \cite{Stein1986} developed exchangeable pairs systematically as a means of constructing Stein identities. I will discuss these ideas of Stein briefly but with a little of my own perspective based on Barbour and Chen \cite{Barbour-Chen2014} and Chen \cite{Chen1998}.

 A pair of random variables $(W,W')$, defined on the same proability space, is said to be an exchangeable pair if $\calL(W,W')=\calL(W',W)$. We discuss two examples of exchangeable pair.

\textit{Example 1.} Let $X_1, \cdots, X_n$ be independent random variables with $\E X_i =0, i = 1,\cdots,n$, and let $W=\sum_{i=1}^n X_i$. To construct $W'$ so that $(W,W')$ forms an exchangeable pair, we let $X_i',\cdots,X_n'$ be an independent copy of $X_1, \cdots, X_n$ and let $I$ be a random integer, uniformly distributed over $\{1,\cdots,n\}$ and independent of all the other random variables. Let $W' = W - X_I + X_I'$. Then $(W,W')$ is an exchangeable pair. 

\textit{Example 2.} This example is a combinatorial sum considered by Stein in the 1960s. Let $n \ge  2$ and let $\left(c_{ij}\right)_{i,j=1}^n$ be an $n\times n$ array of real numbers such that for every $i$, $\sum_{j=1}^nc_{ij} =0$ and for every $j$, $\sum_{i=1}^nc_{ij} = 0$. Let $\pi$ be a random permutation of $(1,\cdots,n)$ and let $W = \sum_{i=1}^n c_{i\pi(i)}$. To define $W'$, let $(I,J)$ be uniformly distributed over $\{(i,j): i\neq j; i,j = 1,\cdots, n\}$ and be independent of $\pi$. Define
\begin{equation} \label{def-pi'}
\pi'= \left\{\begin{array}{rcl}
\pi(i), &\text{if}& i \neq I,J,\\
\pi(J), & \text{if} & i = I,   \\
\pi(I), & \text{if} & i = J.
\end{array}\right.
\end{equation}
Then $\pi'$ is a random permutation of $(1,\cdots,n)$. Define 
$$
W' = \sum_{i=1}^n c_{i\pi'(i)} = W - c_{I\pi(I)}-c_{J\pi(J)}+c_{I\pi(J)}+c_{J\pi(I)}.
$$
Then $(W,W')$ is an exchangeable pair.

Suppose $W$ and $W'$ take values in a set $\calS$ and $(W,W')$ is an exchangeable pair. It follows immediately that for any antisymmetric function $F: \calS^2 \rightarrow \R$, that is, $F(w,w') = - F(w',w)$, such that $\E|F(W,W')| < \infty$, $\E F(W,W') = 0$, and hence also $\E[\E(F(W,W')|W)] =0$. 

In order to convert $\E[\E(F(W,W')|W)] =0$ to $\E Lf(W) = 0$ for some operator $L$ and $f$ in some class of functions, we need to define an appropriate subcollection of antisymmetric functions $F$. To this end, we first observe that if $F:\calS^2 \rightarrow \R$ is antisymmetric, $F(w,w') = F(w,w')/2 - F(w',w)/2$.  So $F$ is antisymmetric if and only if $F(w,w') = \phi(w,w') - \phi(w',w)$ for some $\phi: \calS^2 \rightarrow\R$.  Next we focus on those antisymmetric functions of the form
\begin{equation} \label{antisym-F}
F(w,w') = \psi(w,w')f(w') - \psi(w',w)f(w)
\end{equation}
where $\psi:\calS^2 \rightarrow \R$ and $f: \calS \rightarrow \R$. 

For normal approximation using an exchangeable pair $(W,W')$, Stein \cite{Stein1986} proposed this antisymmetric function
\begin{equation} \label{antisym-F-normal}
F(w,w') = (w'-w)[f(w') + f(w)]
\end{equation}
and assumed
that the exchangeable pair $(W,W')$ satisfies
\begin{equation} \label{regression}
\E(W'|W) = (1-\lambda)W,
\end{equation}
where $0 <\lambda < 1$. Note that the $F$ in (\ref{antisym-F-normal}) corresponds to the case $\psi(w,w') = w' - w$ in (\ref{antisym-F}). Note also that (\ref{regression}) implies $\E W = 0$. Assume $\E W^2 = 1$. Then (\ref{regression}) also implies 
\begin{equation} \label{variance}
\E(W' -W)^2 = 2\lambda\E W^2 = 2\lambda.
\end{equation} 
As $F$ is antisymmetric, for bounded absolutely continuous $f$ with bounded $f'$, 
\begin{equation} \label{EF=0}
\E (W'-W)[f(W') + f(W)] = 0.
\end{equation}
Applying (\ref{regression}) to (\ref{EF=0}), we obtain the Stein identity,
\begin{eqnarray} \label{Stein-id-ex-pair1}
~~~~~\E\{f'(W) - Wf(W)\} &=& \E\left\{\E\left[\left(1 - \frac{(W'-W)^2}{2\lambda}\right)\big|W\right]f'(W)\right\} \\ \notag
&& - \frac{1}{2\lambda}\E(W'-W)^2[f'(W + (W'-W)U_1) - f'(W)], 
\end{eqnarray}
where $U_1$ is a uniformly distributed over $[0,1]$ and independent of $W,W'$.
If we assume $f \in \calC_B^2$, (\ref{Stein-id-ex-pair1}) can be written as
\begin{eqnarray} \label{Stein-id-ex-pair2}
\E\{f'(W) - Wf(W)\} &=& \E\left\{\E\left[\left(1 - \frac{(W'-W)^2}{2\lambda}\right)\big|W\right]f'(W)\right\} \\ \notag
&& - \frac{1}{2\lambda}\E(W'-W)^3U_1f''(W + (W'-W)U_1U_2),
\end{eqnarray} 
where $U_2$ is a uniformly distributed over $[0,1]$ and independent of $W, W', U_1$.
Now let $f = f_h$ in (\ref{Stein-id-ex-pair2}), where $h \in \calC^1$ such that $\|h'\|_\infty \le 1$. Applying Proposition \ref{prop-W} and using $\|f'_h\|_{\infty} \le \|h'\|_{\infty}$ and $\|f''_h\|_{\infty} \le 2\|h'\|_{\infty}$ (see Lemma 2.4 in Chen, Goldstein and Shao \cite{C-G-S2011}), we obtain
\begin{eqnarray} \label{d-W-ex-pair}
d_W(W,Z) &\le& \frac{1}{2\lambda}\sqrt{\Var\{\E\left[(W'-W)^2|W\right]\}} + \frac{1}{2\lambda}\E|W'-W|^3.
\end{eqnarray}
It is less straightforward to obtain a bound of the same order on the Kolmogorov distance. Stein \cite{Stein1986} (Theorem III.1) proved that
$$
 d_K(W,Z) \le \frac{1}{\lambda}\sqrt{\Var\{\E\left[(W'-W)^2|W\right]\}}+ (2\pi)^{-1/4}\sqrt{\frac{1}{\lambda}\E|W'-W|^3}.
$$ 

Let us now consider the antisymmetric function for which $\psi \equiv 1$. In this case, $F(W,W') = f(W') - f(W)$. Assume that $(W,W')$ is an exchangeable pair satisfying (\ref{regression}). Starting with $\E[f(W') - f(W)] =0$ and assuming $f \in \calC_B^3$, we obtain the following Stein identity.
\begin{eqnarray} \label{Stein-id-ex-pair3}
\E\{f''(W) - Wf'(W)\} &=& \E\left\{\E\left[\left(1 - \frac{(W'-W)^2}{2\lambda}\right)\big|W\right]f''(W)\right\} \\ \notag
&& - \frac{1}{\lambda}\E(W'-W)^3U_1^2U_2f'''(W + (W'-W)U_1U_2U_3),
\end{eqnarray}
where $U_1, U_2, U_3$ are independent, uniformly distributed over $[0,1]$, and independent of $W, W'$.
The left hand side of (\ref{Stein-id-ex-pair3}) is of the form $\E \calA f(W)$, where $\calA$ given by $\calA f(w) = f''(w) - wf'(w)$ is the generator of the Ornstein-Uhlenbeck process whose stationary distribution is $\calN(0,1)$. This leads us to the generator approach of Barbour \cite{Barbour1988, Barbour1990}, in which for approximation by the stationary distribution of a Markov process $(\xi(t))_{t \ge 0}$, the Stein equation is $\calA f(w) = h(w) - \E h(\xi)$,
where $\calA$ is the generator of the Markov process, $\xi$ has the stationary distribution, and $h$ belongs to a suitable separating class of functions, $\calH$. A bounded solution $f_h$ of the Stein equation is given by
$$
f_h(w) = - \int_0^\infty \E_w[h(\xi(t)) - \E h(\xi)]dt,
$$
and $d_\calH(W, \xi) = \sup_{h \in \calH}|\E\calA f_h(W)|$. So the Stein identity (\ref{Stein-id-ex-pair3}) can be interpreted as one for the generator approach.

Note that exchangeability is not used in the arguments leading to (\ref{Stein-id-ex-pair3}). We only use $\calL(W) = \calL(W')$ and $\E(W'|W) = (1-\lambda)W$. This fact that a Stein identity using $(W,W')$ can be constructed without exchangeability was observed by R\"ollin \cite{Rollin2008}.

The approach using exchangeable pairs has proved extremely useful in many contexts. In many examples, the exchangeable pair can be realized as a pair of successive states in a stationary reversible Markov chain. In the book Diaconis and Holmes \cite{Diaconis-Holmes2004}, for instance, there are chapters showing how Stein’s method and exchangeable pairs can be effectively exploited in a variety of quite disparate settings. Other examples that have since been influential include those in Rinott and Rotar \cite{Rinott-Rotar1997}. Chen and Fang \cite{Chen-Fang2015-1} used exchangeable pairs and a concentration inequality on a more general version of the combinatorial central theorem and obtained a third moment bound on the Kolmogorov distance. The combination of linear regression condition and the exchangeable pair has been effectively developed in the setting of multivariate normal approximation by Reinert and R\"ollin \cite{Reinert-Rollin2009}.

\section{Epilogue} Recent developments have continued to enlarge and enrich the structure of Stein's method, both in theory and applications, to include distributions such as the Dickman distribution, with applications to sorting algorithms and probabilistic number theory \cite{B-G2019}, that bear less and less resemblance to the classical ones. 

The breadth of its connections with other areas has also widened. Under the framework of sublinear expectation for modeling the uncertainty of probabilities and distributions in real data, Stein's method has been applied to prove a central limit theorem with rate of convergence \cite{F-P-S-S2019}. In applications to the physical sciences, Stein's method has now been used to prove a central limit theorem for the free energy of the random field Ising Model \cite{Chatterjee2019}.

The branch connecting Stein's method to Malliavin calculus has continued to grow and evolve. For example, the control by the fourth moment in Weiner chaos to the quality of normal approximation in that setting is now well understood. The variant taking Poisson spatial processes as input has produced a plethora of tight results in stochastic geometry, such as those for Voronoi tesselations \cite{L-S-Y2019}. In a related offshoot in stochastic analysis, Stein kernels can now be used to obtain improvements in the log Sobolev inequality \cite{L-N-P2015}, and their connections to optimal transport have already begun a fruitful interplay \cite{Fathi2019, Herry2018}.

More recently, in the field of data analysis, the \textit{Kernalized Stein Discrepancy} (KSD) \cite{C-S-G2016, G-M2015, G-M2017, L-L-J2016} and \textit{Stein Variational Gradient Descent} (SVGD) \cite{L-W2016}, are now well known practical machine learning algorithms for data fitting and Monte Carlo type simulation that are based on Stein's original idea of measuring distributional distance using solution of specialized differential equations. In high dimensional data analysis, Stein's ideas are currently applied to extend the advantages of shrinkage estimation to non-Gaussian settings, and to estimate parameters and evaluate the costs for the violation of Gaussian assumptions in single-index models and compressed sensing \cite{F-G-R-S2020}.

After 49 years, Stein's ideas are as vibrant as ever with a momentum of their own. The ever widening and often unexpected branching of the path the method has taken suggests that even though its future direction is unforeseeable, neverthless we are confident that it will continue to lead us to new, productive and exciting areas.

\section*{Acknowledgement} I am thankful to Larry Goldstein and Yu-Kiang Leong for their valuable comments which helped me improve the exposition of this paper. I am particularly thankful to Larry for his contribution to Section 7. Thanks also go to Adrian R\"{o}llin for being always around to render a helping hand whenever I encountered any difficulty with using \LaTeX.

\end{document}